\documentclass[11pt]{article}

\usepackage[margin=1in]{geometry}
\usepackage{amssymb, amsthm}

\usepackage{amsmath, bm}
\usepackage{graphicx}
\usepackage{float}
\usepackage{soul, color}
\usepackage[font=scriptsize, labelfont=bf]{caption}
\usepackage{enumitem}

\usepackage[round]{natbib}
\newcommand{\argmax}{\operatornamewithlimits{argmax}}

\newcommand{\rev}[1]{{\color{black}{#1}}}
\newcommand{\x}{\mathbf{x}}
\renewcommand{\v}{\mathbf{v}}
\newcommand{\mc}{\mathcal}

\title{Markov Decision Process Design: A Framework for \\ Integrating Strategic and Operational Decisions}
\author{Seth Brown\thanks{Computational Applied Mathematics \& Operations Research, Rice University, Houston, TX, USA}
\and Saumya Sinha\thanks{Industrial \& Systems Engineering, University of Minnesota, Minneapolis, MN, USA (Corresponding author)} 
\and Andrew Schaefer\footnotemark[1]}
\date{}

\begin{document}
\maketitle

\begin{abstract}
We consider the problem of optimally designing a system for repeated use under uncertainty. We develop a modeling framework that integrates design and operational phases, which are represented by a mixed-integer program and discounted-cost infinite-horizon Markov decision processes, respectively. We seek to simultaneously minimize the design costs and the subsequent expected operational costs. 
This problem setting arises naturally in several application areas, as we illustrate through examples.
We derive a bilevel mixed-integer linear programming formulation for the problem and perform a computational study to demonstrate that realistic instances can be solved numerically.
\end{abstract}

\bigskip

\section{Introduction}

We present a new modeling framework for the problem of optimally designing a system for repeated use under uncertainty. Consider an agent who seeks to make decisions in two phases, namely, a design phase followed by an operational phase. The design phase corresponds to one-time strategic decisions that have long implementation cycles and are difficult to modify once executed. As such, these are here-and-now static decisions that are made upfront. The operational phase begins once design decisions have been executed and any uncertainty in their outcomes has been realized. This phase corresponds to day-to-day dynamic decisions that are routinely made in an uncertain environment once the system is in use. 

Examples of this setting arise naturally in several application areas. For instance, consider a firm that seeks to acquire a new fleet of vehicles. In the design phase, the firm determines the number and specifications of the vehicles to purchase, subject to budgetary constraints. The purchase orders may take years to fulfill and are difficult to revise. Moreover, the actual performance of the delivered vehicles is uncertain and unknown at the time of procurement. The operational phase corresponds to the routine deployment of the fleet after vehicle delivery. As another example, consider a retailer who seeks to open new retail outlets at a few locations selected from a set of available choices. The design phase corresponds to decisions about the number, location, and specifications (e.g., size and capacity) of the new outlets, while the operational phase corresponds to inventory management at the retail outlets once they are ready for use. 

We highlight two key features of this decision problem: 
\begin{enumerate}[label=(\roman*), topsep=0pt, itemsep=-2pt]
\item interdependence of the decision phases: As evident from the examples, the two phases of decision-making are naturally interdependent -- optimal design decisions must factor in the long-term operational cost of the enterprise, while the parameters of the operational problems are determined by design phase decisions. 
\item two sources of uncertainty: As first noted by \citet{bs06}, two distinct sources of uncertainty are present in this decision framework. Design uncertainty refers to uncertainty in the outcome of the design phase (e.g., variability in vehicle performance after purchase), and is not realized until the design phase has been executed. Operational uncertainty refers to the stochastic environment in which operational decisions are made (e.g., uncertain demand in the inventory management problem). 
\end{enumerate}

Despite the natural interdependence, standard modeling approaches typically view the two phases as disparate sub-components, due in part to the complexity of their relationship and high dimensionality of the resulting problem \cite{burnak2019integrated}.
Several authors have discussed the need to integrate the different phases of decision-making \cite{burnak2019integrated, grossmann2005enterprise, avraamidou2019bi, badejo2022integrating}, noting that synergistic interactions across different phases can lead to gains in efficiency and profits. In this work, we present an optimization model that captures the hierarchy of these decisions, and incorporates both types of uncertainty. 
We model the static design phase as a mixed-integer program (MIP), and the dynamic operational decisions via infinite-horizon discounted cost Markov decision processes (MDPs).
Design uncertainty is modeled via a finite set of scenarios, and operational uncertainty through probabilistic state transitions in the MDPs. Then, we seek to minimize the combined expected cost of the MIP and the subsequent discounted-cost MDPs that are determined after the first-stage decisions have been made. 
This problem may also be viewed as one where we choose optimal design parameters for MDPs so as to minimize the expected long-run operational cost under a family of uncertain scenarios. 
To the best of our knowledge, we are the first to consider the problem of optimal MDP design.

\citet{bs06} consider an {\it adversarial} version of this problem, where an agent first makes a design decision to minimize costs, and an opponent subsequently makes operational decisions to maximize the long-run costs to induce the greatest disruption to the agent. The contradiction of goals in the two stages allows for the reformulation of the problem as a stochastic MIP with continuous second-stage variables. 
In contrast, we follow a {\it cooperative} approach where the same agent makes the decisions in both stages, and the objectives in both stages are aligned so that we simultaneously minimize both the design and operational costs.
This alignment renders the stochastic MIP formulation of \cite{bs06} no longer achievable. Instead, we formulate the problem as a bilevel program.
Note that this contrast between the cooperative and adversarial approaches is different from the usual distinction between `optimistic' and `pessimistic' formulations in bilevel optimization \cite{dempe2015bilevel}. In our framework, the decision-maker in the strategic phase (i.e., the leader problem) is only interested in the optimal values of the follower problems, and is unaffected by the particular choice of the followers' optimal solutions. As such, questions about optimistic and pessimistic formulations do not arise. 

Bilevel optimization has been used in several domains to integrate decisions that were previously modeled independently,
such as facility location and routing \cite{gonzalez2022routing, tordecilla2023routing}, and design and scheduling of production systems \cite{avraamidou2019bi, leenders2023bilevel}. To our knowledge, we are the first to consider the general problem of integrating the design and operational phases where the latter is modeled as an MDP, incorporating both types of uncertainty.

Recent years have witnessed significant advances in the development of solution approaches for bilevel programming problems, which are generally characterized by their (non-)linearity and the presence of discrete and/or continuous variables. 
Our formulation consists of a linear bilevel program with a mixed-integer leader problem and multiple follower problems with continuous variables.
\citet{ah97} show that linear bilevel programs can be reformulated as mixed binary linear programs, and 
\citet{fl16} describe a general approach that focuses on adapting existing MIP solvers to tackle bilevel problems. 
\citet{si09} present a Benders-like decomposition approach for linear problems, 
and \citet{xw14} solve the problem with bounded integer leader variables.
\citet{ag88} discusses a decentralized bilevel problem that has a single leader and multiple followers, which is one way of viewing the stochastic variant of a bilevel programming problem. 
\citet{tr16} describe the implementation of the bilevel solution software MibS, which uses a branch-and-cut approach to solve mixed-integer linear bilevel programs. 
Several strategies for solving fully general bilevel MIPs have been developed, many of which are discussed by \citet{fl17} in the context of their software.
\citet{sc13} survey solution approaches for linear bilevel programs. In particular, the authors discuss three categories of methods for mixed-integer problems: reformulation, branch-and-bound/branch-and-cut, and parametric programming. 
\citet{kleinert2021survey} provide a more recent survey of solution algorithms, including a discussion of some nonlinear cases.

\label{rev:intro}\rev{We note that the focus of this paper is to present a general modeling framework that (a) applies to a large class of problems that arise in many application domains, and (b) is simple enough to be solved directly using off-the-shelf bilevel optimization solvers. Specialized solution algorithms that leverage application-specific problem structure may be more efficient in some instances, but are beyond the scope of the present work.}
In Section \ref{sec:model}, we describe the problem setup in detail, and derive a bilevel linear MIP formulation. 
In Section \ref{sec:examples}, we provide examples to illustrate the broad applicability of the proposed modeling environment. In Section \ref{sec:numerics}, we present a numerical study to demonstrate the feasibility of this framework for realistically-sized problems using existing optimization software.

\section{MDP Design: An Integrated Framework for Design and Operations}
\label{sec:model}
We consider a two-phase integrated decision framework, consisting of a design phase and an operational phase. In the static design phase, an agent chooses optimal design variables modeled via a MIP, and during the operational phase, she determines optimal decision rules for an infinite-horizon MDP that depends on the design decisions. We use a scenario approach to account for uncertainty in design execution between the two phases, whereas operational uncertainty is naturally modeled via probabilistic state transitions in the MDPs. 
Our objective is to minimize the combined expected cost of the design-phase MIP and the subsequent discounted-cost MDPs that represent the operational phase, taking into account both types of uncertainty. Dependence of the MDP on design decisions is modeled through the MDP cost functions, which are assumed to be affine functions of the leader variables. 
This dependence may be extended to settings where the MDP state and action spaces depend on the design decisions by means of a big-M approach that preserves the linear relationship; see Section \ref{sec:inventory} for an example. 

\paragraph{Design Problem}
Let $\x \in \mathcal{X} \subseteq \mathbb{R}^{n_1} \times \mathbb{Z}^{n_2}$ represent the mixed-integer design variables. The design MIP is
parameterized by the cost vector $\mathbf{c} \in \mathbb{R}^{n}$ (where $n = n_1+n_2$), constraint matrix $\mathbf{A} \in \mathbb{R}^{m \times n}$, and constraint right-hand-side $\mathbf{b} \in \mathbb{R}^{m}$.
To model uncertainty in the outcomes of the design phase, we consider a finite set of scenarios $\mathcal{K}$, where each scenario $k \in \mathcal{K}$ occurs with probability $q_k > 0$. 
Thus, $\mathbf{q}\in [0,1]^{\vert \mathcal{K}\vert}$ and $\sum\limits_{k \in \mathcal{K}} q_k = 1$. Note that we use boldface letters (e.g., $\mathbf{x}$) to denote vectors, and subscripted regular typeface letters (e.g., $x_i$) to denote their components; matrices are denoted by boldface upper case letters.
Let $u^k(\x)$ be the optimal expected total cost of the operational phase in scenario $k$ for a given design decision $\x$. Then, the design problem is expressed as 
\begin{equation}
\label{eqn:design}
    \tag{MIP}
\begin{aligned}
    \min & ~ \mathbf{c}^\top \x + \sum\limits_{k \in \mathcal{K}} q^k u^k(\x) \\
    \text{s.t. } & ~ \mathbf{A} \x = \mathbf{b}, \\
    & ~ \x \in \mathcal{X}.
\end{aligned}
\end{equation}

\paragraph{Operational Problem}
For a fixed design decision $\x$ and scenario $k \in \mathcal{K}$, we model the operational phase as an infinite-horizon discounted cost MDP with a (possibly scenario-dependent) discount factor $\lambda_k \in (0,1)$, finite state space $\mathcal{S}^k$, and finite action space $\mathcal{A}^k$. For each state-action pair $(s,a) \in \mathcal{S}^k \times \mathcal{A}^k$, the agent incurs an immediate cost $(\mathbf{f}_{s,a}^k)^\top\x + g_{s,a}^k$ that depends linearly on $\x$. Transition probabilities are denoted by $\mathbf{p}^k( \cdot \,\vert \, s, a)$. 
A (stationary deterministic) decision rule $\mathbf{d}^k $ is a function that assigns one action to every state in $\mathcal{S}^k$.
Given a decision rule $\mathbf{d}^k$, let $\mathbf{P}_{\mathbf{d}^k}$ be the associated transition probability matrix and $\mathbf{h}_{\mathbf{d}^k}(\x)$ be the cost vector whose $s^{\textrm{th}}$ entry is the cost for state-action pair $(s, d^k_s)$. Then, the expected value of policy $\mathbf{d}^k$ is
\begin{equation}
\label{eqn:val_policy}
\v_{\mathbf{d}^k}(\x) = (\mathbf{I}-\lambda_k\mathbf{P}_{\mathbf{d}^k})^{-1}\mathbf{h}_{\mathbf{d}^k}(\x),
\end{equation}
where $\mathbf{I}$ denotes the identity matrix.
The optimal value function of the MDP, denoted by $\v^k(\x)$, is obtained by minimizing \eqref{eqn:val_policy} componentwise over all possible decision-rules $\mathbf{d}^k$ (or equivalently, given positive weights $\bm{\beta}$, by minimizing $\bm{\beta}^\top\v_{\mathbf{d}^k}$ \cite{pm14}). 


Standard methods for solving infinite-horizon MDPs include value iteration, policy iteration, and linear programming (LP) \citep{pm14}. While the first two are more popular for solving stand-alone MDPs, we use the LP approach to formulate the integrated problem as a bilevel program whose second level is an LP. 
%
The LP method uses Bellman's equations to derive linear constraints on the optimal values of states, so that $\v^k(\x)$ is the optimal solution of \eqref{eqn:LP}.
\begin{equation}
\label{eqn:LP}\tag{LP($\x,k$)}
\begin{aligned}
\max_{\v\in \mathbb{R}^{\vert \mathcal{S}^k \vert}} & ~ \mathbf{1}^\top\v \\
\text{s.t. } & ~ v_s \leq (\mathbf{f}_{s,a}^k)^\top\x+g_{s,a}^k + \lambda_k\sum_{s'\in \mathcal{S}^k} p^k(s' \vert\, s,a)\ v_{s'}, \quad \forall\ a \in \mathcal{A}^k, s \in \mathcal{S}^k,
\end{aligned}
\end{equation}
where $\mathbf{1}$ is the vector of all ones.
Let $\alpha^k_s$ be the probability that the MDP in scenario $k$ initially occupies state $s \in \mathcal{S}^k$. Then, in the notation of \eqref{eqn:design}, we have $u^k = (\bm{\alpha}^k)^\top \v^k(\x)$ for all $k \in \mc{K}$. 

The objective function of \eqref{eqn:design} is the sum of design costs and the expected {\it optimal} costs from the MDPs. Unlike the problem with adversarial recourse in \cite{bs06}, the function $\v^k(\x)$ is not convex; instead, it is piecewise linear concave \cite{bj91}. 
To see this, note from Equation \eqref{eqn:val_policy} that $\v^k(\x)$ is the minimum of a family of linear functions of $\x$ and therefore piecewise linear and concave \cite{boyd2004convex}. In fact, it continues to be concave in the more general case where immediate costs of the MDP are concave (and not necessarily affine) functions of the design variables $\x$.

A brute-force approach to solving \eqref{eqn:design} consists of enumerating all its feasible solutions and solving $|\mathcal{K}|$ MDPs to optimality for each objective function evaluation. 
As such, this approach may be feasible when $\mathcal{X}$ is a small finite set, but becomes impractical when it is hard/expensive to enumerate all the feasible solutions to \eqref{eqn:design}, or if the design phase includes continuous variables.

\subsection{Bilevel Programming Formulation}
\label{sec:bilevel}

To derive tractable reformulations of \eqref{eqn:design}, we replace the functions $\v^k(\x)$ with auxiliary variables $\v^k$, and explicitly constrain these variables to equal the optimal values of the corresponding MDPs. 
Thus, we model \eqref{eqn:design} as a (mixed-integer) linear bilevel program with a single leader corresponding to the choice of the design variable $\x$, and $| \mathcal{K}|$ follower problems, each comprised of the LP for the associated MDP. As such, the decision variables for the $k^{\textrm{th}}$ follower problem are $v_s^k$, the value of state $s$ in scenario $k$, for all states $s \in \mathcal{S}^k$. 
Therefore, a bilevel programming formulation for \eqref{eqn:design} is given by 
\begin{equation}
\label{eqn:bilevel}
\begin{aligned}
\min_{\x, \v^k} &~ \mathbf{c}^\top\x + \sum_{k \in \mathcal{K}}q^k [(\bm{\alpha}^k)^\top \v^k]\\
\text{s.t. } &~ \mathbf{A}\x = \mathbf{b},\\ 
&~ \x \in \mathcal{X}, \\
&~ \v^k \in \argmax\limits_{\v \in \mathbb{R}^{\vert \mathcal{S}^k \vert}} \Big\{ \mathbf{1}^\top\v : v_s \leq (\mathbf{f}_{s,a}^k)^\top\x+g_{s,a}^k + 
\lambda_k\sum_{s'\in \mathcal{S}^k} p^k(s'\vert s,a)v_{s'}, \ \forall\, a \in \mathcal{A}^k, s \in \mathcal{S}^k \Big\}, \quad \forall\, k \in \mathcal{K}.
\end{aligned}
\end{equation}

Problem \eqref{eqn:bilevel} is a linear bilevel program with a mixed-integer leader problem and $|\mathcal{K}|$ independent continuous follower problems---an equivalent formulation with a single follower is obtained by moving the expectation from the leader's objective to the follower's objective and using a single follower variable $v = \sum_{k \in \mathcal{K}}q^k ((\bm{\alpha}^k)^\top \v^k)$. 
Because the MDPs have unique optimal value functions \cite{pm14}, the question of optimistic versus pessimistic formulations of the bilevel problem (e.g., see \citep{cm05}) do not arise even though the followers' optimal policies may not be unique. 

\label{rev:bilevelnphard}
\rev{Linear bilevel optimization is known to be strongly NP-hard \cite{hansen1992new}}. In recent years, several researchers have proposed algorithms and developed open-source software for solving general mixed-integer linear bilevel programs; see \cite{sc13, kleinert2021survey} for an overview. As such, \eqref{eqn:bilevel} can be solved using existing bilevel solvers. 
\subsection{MIP Formulation}
To obviate the need for specialized software or algorithms for bilevel programs, a common solution approach utilizes the dual of the follower problem to reformulate \eqref{eqn:bilevel} as a single-level MIP. The dual of \eqref{eqn:LP} is given by
\begin{equation}
\label{eqn:LPdual}
\begin{aligned}
\min &~ \sum_{s\in \mathcal{S}^k} \sum_{a\in \mathcal{A}^k} \gamma_{s,a}^k \cdot
((\mathbf{f}_{s,a}^k)^\top\x + g_{s,a}^k)\\
\text{s.t. } &~ \sum_{a \in \mathcal{A}^k} \Big[\gamma_{s,a}^k - \lambda \sum_{s' \in \mathcal{S}^k} \gamma_{s',a}^k p^k(s\,\vert\, s', a) \Big] = 1 \quad \forall \ s \in \mathcal{S}^k, \\
&~ \gamma_{s,a}^k \geq 0 \quad \forall \ s \in \mathcal{S}^k, a \in \mathcal{A}^k.
\end{aligned}
\end{equation}

Then, we derive a nonlinear MIP reformulation of \eqref{eqn:bilevel} by substituting the dual problems \eqref{eqn:LPdual} for the follower problems, and including complementary slackness conditions as (nonlinear) constraints. 
Moreover, we can linearize the MIP using a big-M approach as follows \cite{fortuny1981representation}:
\begin{equation}
\label{eqn:mip_reform}
\begin{aligned}
\min\ & \mathbf{c}^\top\x+\sum_{k\in \mathcal{K}}q^k\sum_{s\in \mathcal{S}^k}\alpha_s^kv^k_s\\
\text{s.t. } & \mathbf{A}\x = \mathbf{b},\\
& \x \in \mathcal{X},\\
& v_s^k \leq (\mathbf{f}_{s,a}^k)^\top\x + g_{s,a}^k+\lambda_k\sum_{j\in \mathcal{S}}p^k(j\vert s,a)v_j^k \quad \forall s \in \mathcal{S}^k, a \in \mathcal{A}^k, k \in \mathcal{K},\\
& \sum_{a \in \mathcal{A}^k} [\gamma_{s,a,k} -\lambda_k\sum_{j \in \mathcal{S}^k} \gamma_{a,j,k}p^k(s\,\vert\,j, a)] = q^k\alpha_s^k \quad \forall \ s \in \mathcal{S}^k,\ k \in \mathcal{K},\\
& \gamma_{s,a,k} \leq M_k \delta_{s,a,k} \quad \forall \ s \in \mathcal{S}^k, \ a \in \mathcal{A}^k, \ k \in \mathcal{K},\\
& (\mathbf{f}_{s,a}^k)^\top\x + g_{s,a}^k-v_s^k+\lambda_k\sum_{j\in \mathcal{S}}p^k(j\vert s,a)v_j^k \leq M_k' (1-\delta_{s,a,k})  \quad \forall \ s \in \mathcal{S}^k, \ a \in \mathcal{A}^k, \ k \in \mathcal{K},\\
& \gamma_{s,a,k}\geq 0, \ \delta_{s,a,k} \in \{0,1\} \quad \forall \ s \in \mathcal{S}^k, \ a \in \mathcal{A}^k, \ k \in \mathcal{K}.
\end{aligned}
\end{equation}

\rev{
This single-level MIP reformulation is fairly common in the bilevel optimization literature.}
It allows us to solve the integrated problem \eqref{eqn:bilevel} using readily available off-the-shelf MIP solvers. 
Moreover, as shown in Section \ref{sec:inventory}, this approach can provide additional modeling flexibility in cases where the state and action spaces for the follower depend on the design decision. 
\rev{However, the computational performance of this approach depends on the selection of suitable big-M parameters, which is a non-trivial question. \citet{kleinert2020nofreelunch} show that there is no polynomial time general-purpose method for selecting a correct big-M unless P$=$NP. 
Furthermore, \citet{kleinert2023there} provide examples where a poor choice of big-M parameters can lead to arbitrarily bad solutions. Even so, many authors have used problem context to carefully select parameters that perform well in practice; see \cite{sc13, kleinert2021survey, kleinert2023there} and references therein.
}\label{rev:bigMcomplexity}

In formulation \eqref{eqn:mip_reform}, for constraints corresponding to state $s\in \mc{S}^k$, we can choose $M_k$ to equal $M_{s,k}$, the maximum number of (discounted) times state $s$ could be reached in scenario $k$. Similarly, we can set $M_k'$ to equal $M'_{s,k}$, the most costly possible outcome starting in state $s$ in scenario $k$ for any feasible $\x$. The latter is found by solving an LP for each state-scenario pair. More simply, we may choose $M_{k} =1/(1-\lambda_k)$ and $M'_{k} = t^k/(1-\lambda_k)$, where $t^k$ is the maximum cost for any state-action pair in scenario $k$. \label{rev:bigm}
%
However, with this choice of big-M parameters, our preliminary computational tests indicated poorer computational performance for the MIP approach than for directly solving the bilevel model using specialized bilevel optimization software. 
\rev{Further sharpening of these parameters may be possible in specific application contexts, but is beyond the scope of the present work.}
Therefore, we focus on the bilevel model in our numerical experiments in Section \ref{sec:numerics}.

\label{rev:extensions}
\rev{Finally, we note that the above framework can be extended to more general model settings. 
Specifically, we can directly apply the approach of Section \ref{sec:bilevel} to other forms of decision-models for the operational phase, provided they have a suitable linear programming representation, such as finite-horizon total cost MDPs \cite{bhattacharya2017linear} and finite-horizon risk-sensitive MDPs \cite{kumar2015finite}. However, because the lower level LP now has a decision-variable $v_{s,t}$ for each state $s$ and period $t$, the resulting follower problems will be larger in size, scaling linearly with the length of the horizon. 
From a modeling standpoint, we can also extend the framework to cases where the MDP cost function depends nonlinearly on the design variables or where transition probabilities depend on the design decisions, but the resulting bilevel optimization problems will no longer be linear. As such, our ability to solve these problems may be limited. }

\section{Applications}
\label{sec:examples}
We now illustrate that this problem setting arises in several application areas with interdependent design and operational phases.

\subsection{Reliability}
In designing a system for long-term operation, minimizing downtime and maintenance costs is an important consideration which must be balanced against the need to lower initial costs. As such, there is a tradeoff between, for example, selecting expensive high-quality components to build the system, or using low-cost components that lead to increased system downtime and maintenance requirements during regular operation. An integrated model would make optimal design choices to balance this tradeoff. 

\citet{tc95} consider this problem and propose two solution approaches. The first aims to iteratively solve the MDP for a range of design decisions, which cannot be guaranteed to recover the true optimal unless the number of design choices is finite. The second reformulates the problem as a single nonlinear non-convex problem without structure that would lend itself to tractability.  
We discuss this problem in the context of our framework. Specifically, we use MDPs to model system reliability and extend the standard approach to additionally consider the problem of purchasing the best components from a set of options under budgetary constraints.

Let $T$ be a system dependent on some subset of $N$ components arranged in series, and let $B$ be a boolean expression for whether $T$ is operational based on the current present and working components. This can easily be expanded to a network of agents, each of which has its own components and may or may not be functioning at a given time, with an arbitrary degree of dependence on the other agents to determine the transition probabilities; for the purposes of this discussion, we use a single-system model. 

We model the design problem via mixed-integer decisions $\x \in \mathbb{R}^{\vert N\vert}$ for the amount of each component to purchase (where a continuous decision variable may indicate a volume or mass of some component) at costs $\mathbf{c}$, such that system operation is possible and the purchase cost does not exceed a budget $b$. Let $\mathcal{X}$ be the set of component choices that allow for system operation.
The uncertainty in design outcomes represents the variability in the quality of the components that are eventually delivered, as well as external sources of uncertainty, such as demand for the service provided and the labor market for performing repairs.

In the subsequent operational stage, we define status variables $\mathbf{y} \in \mathbb{R}^{\vert N\vert}$ to indicate the state of each component selected in the first stage, and the system state $s$ is the product of these status variables. Transition probabilities $p_i^k(y_i \,\vert\, s, a)$ for each component $i$ depend on both the current state and the design uncertainty. In any epoch, an action $a$ corresponds a to subset of the components to repair, and the immediate cost $(\mathbf{f}_{s,a}^k)^\top\x + g_{s,a}^k$ includes both the cost of repair and system operation.  

We assume that the system is initially in the fully functioning state denoted by $s_0$. 
Then the bilevel problem takes the form:
\begin{align*}
\min \ & \mathbf{c}^\top\x+\sum_{k \in \mathcal{K}}q^k v_{s_0}(\x,k) \\
\text{s.t. } & \mathbf{c}^\top\x \leq b, \\
& \x \in \mathcal{X}, \\
& \v(\x, k) \in \argmax\ \eqref{eqn:LP}.
\end{align*}


\subsection{Inventory management}
\label{sec:inventory}
Consider an extended version of the inventory management problem, where the first phase corresponds to business startup location decisions in a particular market, and the second phase is a classic inventory management problem with Markovian demand (see \cite{wh93} for a survey). We assume that the costs of ordering and holding inventory are location-dependent. The discount factor accounts for rental or other location dependent costs that are fixed for each location decision and are incurred at each time step.

For the design phase, suppose the agent needs to select $r$ locations out of $R$ available choices. Let $\x \in \mathbb{B}^R$ correspond to binary decisions indicating whether each location was selected, with associated startup costs $\mathbf{c}$. Let $\mathbf{m}$ be the inventory capacities for all locations. We assume that each location can initially order up to a maximum of their capacity. Let $\mathbf{u}$ be the initial inventory ordered, with per unit costs $\mathbf{b}$. Then the objective for the design-only problem would be
\[
\min_{\x \in \mathbb{B}^{R},\, \mathbf{u}\in \mathbb{R}^{R}}~ \mathbf{c}^\top\x+\mathbf{b}^\top\mathbf{u},
\]
with constraints 
\[ 
\mathbf{1}^\top\x = r \quad \text{and } \quad 0 \leq u_i \leq m_ix_i \quad \forall \ i \in \{1, \ldots, R\}.
\]

Design uncertainty may reflect variability in the economic environment and market demand by the time of completion. Once the locations are ready for use, the operational phase corresponds to inventory management at each location, which is modeled as an MDP; see \cite[Section 3.2]{pm14} for an example. Let $\mathbf{f}^k$ be the vector of fixed costs incurred at each location in every period. Given an inventory level $s$ and order quantity $a$, let $\mathbf{o}_{s,a}^k(\x)$ be the ordering costs, $\mathbf{h}_{s,a}^k(\x)$ be the holding costs, and $\mathbf{r}_{s,a}^k(\x)$ be the total sale revenues. We assume that demand is not backlogged, but lost demand incurs shortage costs $\mathbf{p}_{s,a}^k(\x)$. We further assume that the costs and revenues are linear in $\x$. To see why this assumption is justified, fix $i \in \{1, \ldots, R \}$ and note that $x_i$ indicates if the location is operational. Then, if the per unit ordering cost is $o_i^k$, the total ordering cost is given by
\[
o_i^k x_i a_i + M_i^k(1-x_i)a_i = (o_i^ka_i-M_i^ka_i)x_i + M_i^ka_i,
\]
where $M_i$ is a penalty parameter that enforces $a_i=0$ if the $i$-th location has not been selected in the design phase. Thus, the ordering cost is linear in the leader variable. Similar arguments apply to the other cost functions as well. However, there is an added subtlety here because the state and action spaces depend on $\x$. Specifically, the inventory at the $i$-th location in any period must be less than $m_ix_i$, and the order quantities in state $s$ may not exceed $m_ix_i{-}s$. We present two ways of modeling these constraints within the bilevel framework.
\begin{enumerate}[leftmargin=*]
\item Penalty term in cost function: The immediate cost function of the MDP may include additional terms of the form
\[
M_1(m_ix_i-s_i)^+ \text{ and } M_2(m_ix_i-s-a)^+,
\]
where $M_1$ and $M_2$ are (large) penalty parameters. In this case, however, the cost functions will be piecewise linear and convex in $\x$, which makes it more challenging to solve the bilevel problem.

\item Auxiliary variables in the dual problem \eqref{eqn:LPdual}: This is a similar idea as above but utilizes the dual formulation for the MDP. We define auxiliary variables $\rho_s = (m_ix_i-s_i)^+$ and $\mu_{s,a} = (m_ix_i-s-a)^+$ and add terms $M_1 \rho_s$ and $M_2 \mu_{s,a}$ to the objective function of \eqref{eqn:LPdual}. The definitions of $\rho_s$ and $\mu_{s,a}$ may be enforced via linear constraints. This results in a larger dual problem that continues to be linear. Thus, this is an instance where reformulating the bilevel problem as a single-level MIP may be advantageous. 

\end{enumerate}

Dependence of the operational phase on initial inventory $\mathbf{u}$ may be modeled similarly. Thus, our framework may be applied here.

\subsection{Queue Design and Control}

Our third example considers the problem of optimally designing a queueing system. Several authors have studied the optimal selection of queue parameters, such as arrival rates and service times, where performance measures of the system are expressible as functions of these parameters in closed form; see \cite{tc05} for a survey. Our framework allows for greater flexibility by explicitly modeling the operational phase as an MDP. 

In the design phase, an agent seeks to determine the number and types of servers to be installed. Given $n$ server types, let $\x\in \mathbb{Z}_n^+$ denote the number of servers of each type that are selected. Let $\mathbf{r}$ be the maximum allowable number of each type of server, and $t$ be the maximum total number of servers. Then, if the one-time recruitment and training costs per server are given by $\mathbf{c}$, the design problem takes the form
\begin{align*}
\min\limits_{\x \in \mathbb{Z}^n_+} \ & \mathbf{c}^\top\x \\
\text{s.t. } & \x \leq \mathbf{r} \\
& \mathbf{1}^\top\x \leq t.
\end{align*}

The integrated framework additionally considers the uncertainty in the outcome of the design phase and the expected long-run operational cost. 
Design uncertainty in this setting may correspond to the realized service rates of each server type, as well as the servers’ ability to successfully resolve the needs of various customer types. In the operational phase, $J$ types of customers arrive according to some time-discretized arrival process with rate $\eta_j$, $j =1 , \ldots, J$. The state $s$ of the MDP represents the number of customers of each type that are present in the queue. We assume that the system has a finite capacity, so that the state space is finite. However, a penalty $p$ is incurred when the system is at capacity and customers are turned away. 

In every decision epoch, the action $a$ determines the number of customers of each type that are assigned to each server-type. Let $p_{ij}$ be the probability with which a server of type $i$ successfully serves a customer of type $j$, earning a reward $r_{ij}$. We assume that customers exit the system after service regardless of whether service was successful. Let $o_i$ be the operating cost per period of a type-$i$ server. As in Section \ref{sec:inventory}, we may express the costs and rewards as linear functions of the design variables. Then, the agent seeks to determine an assignment strategy that minimizes the total discounted cost less the expected reward. As such, the long-run cost of operating the queueing system can be included in the design decisions using our bilevel framework.

\section{Numerical Results}
\label{sec:numerics}

\begin{table}
\centering

\begin{tabular}{|c|c|c|c|c|r|}\hline
$n$ & $m$ & $|\mathcal{K}|$ & $|\mathcal{S}|$ & $|\mathcal{A}|$ & Avg. solve  time (s) \\ \hline
{\bf 20} & 40 & 20 & 10 & 20 & 3.384\\
{\bf 40} &	40&	20&	10&	20&	5.060\\
{\bf 80}&	40&	20&	10&	20&	7.732\\
{\bf 160} &	40&	20&	10&	20&	12.584\\
{\bf 320} &	40&	20&	10&	20&	34.005\\ \hline
80&	{\bf 10} &	20&	10&	20&	7.431\\
80&	{\bf 20} &	20&	10&	20&	7.464\\
80&	{\bf 40} &	20&	10&	20&	7.439\\
80&	{\bf 80}&	20&	10&	20&	7.640\\
80&	{\bf 160}&	20&	10&	20&	7.664\\ \hline
80&	40&	{\bf 5}&	10&	20&	1.323\\
80&	40&	{\bf 10}&	10&	20&	2.849\\
80&	40&	{\bf 20}&	10&	20&	7.239\\
80&	40&	{\bf 40}&	10&	20&	18.315\\
80&	40&	{\bf 80}&	10&	20&	50.925\\ \hline
80&	40&	20&	{\bf 2}&	20&	0.572\\
80&	40&	20&	{\bf 4}&	20&	1.453\\
80&	40&	20&	{\bf 8}&	20&	4.620\\
80&	40&	20&	{\bf 16}&	20&	25.932\\
80&	40&	20&	{\bf 32}&	20&	237.272\\ \hline
80&	40&	20&	10&	{\bf 5}&	0.772\\
80&	40&	20&	10&	{\bf 10}&	2.524\\
80&	40&	20&	10&	{\bf 20}&	8.005\\
80&	40&	20&	10&	{\bf 40}&	29.119\\
80&	40&	20&	10&	{\bf 80} &	95.440\\ \hline
\end{tabular}
\caption{Average solve times (in seconds) over 100 instances for 25 different combinations of parameter values. The $n$ leader variables are evenly split between binary variables and general integer variables.}
\label{table:solvetimes}
\end{table}

In this section, we illustrate how existing bilevel optimization algorithms may be used to solve realistic instances of the bilevel programs that arise in our framework.
We present results for a general form of the problem, analyzing the performance of the bilevel MIP solver of \citet{fl17} (available at \texttt{https://msinnl.github.io/pages/bilevel.html}) to solve random instances of \eqref{eqn:bilevel} for a range of parameter values. For each set of parameter values in Table \ref{table:solvetimes}, we generated and solved 100 test problems. The cost and constraint coefficients and other parameters for these problems were randomly generated given fixed mean values and variances that are listed in Table \ref{table:parameters}. Relative probabilities were randomly generated as noted in Table \ref{table:parameters} and then normalized to sum to 1. We assumed that all leader variables are integer or binary and the number of integer and binary variables remained the same in each case.

\begin{table}[t]
    \centering
    \begin{small}
    \begin{tabular}{|p{0.45\textwidth}|p{0.15\textwidth}|p{0.15\textwidth}|}\hline
        Parameter & Mean & Variance \\ \hline
        Leader integer variable upper bound & 12 & 1 \\ \hline
        Leader constraint coefficients & $\bar{a} = 12$ & $\sigma = 0.5$ \\ \hline
        Leader constraint upper bound & $\displaystyle \bar{b} = \frac{ (n_1+n_2)}{\bar{a}+2\sigma}$ & $\bar{b}/6$ \\ \hline
        Scenario relative probability & 1 & 0.2 \\ \hline
        Initial state relative probability & 1 & 0.2 \\ \hline
        Transition relative probability & 1 & 0.4 \\ \hline
    \end{tabular}
    \end{small}
    \caption{Parameters of the normal distributions used for instance generation in the numerical experiments. Discount factors were sampled uniformly at random over the interval $[0.92, 0.97]$, leader objective coefficients over $[10,100]$, cost function coefficients over $[-1,1]$, and cost function constant terms over $[10, 40]$.}
    \label{table:parameters}
\end{table}

\begin{figure}
    \centering
    \includegraphics[width=0.55\textwidth]{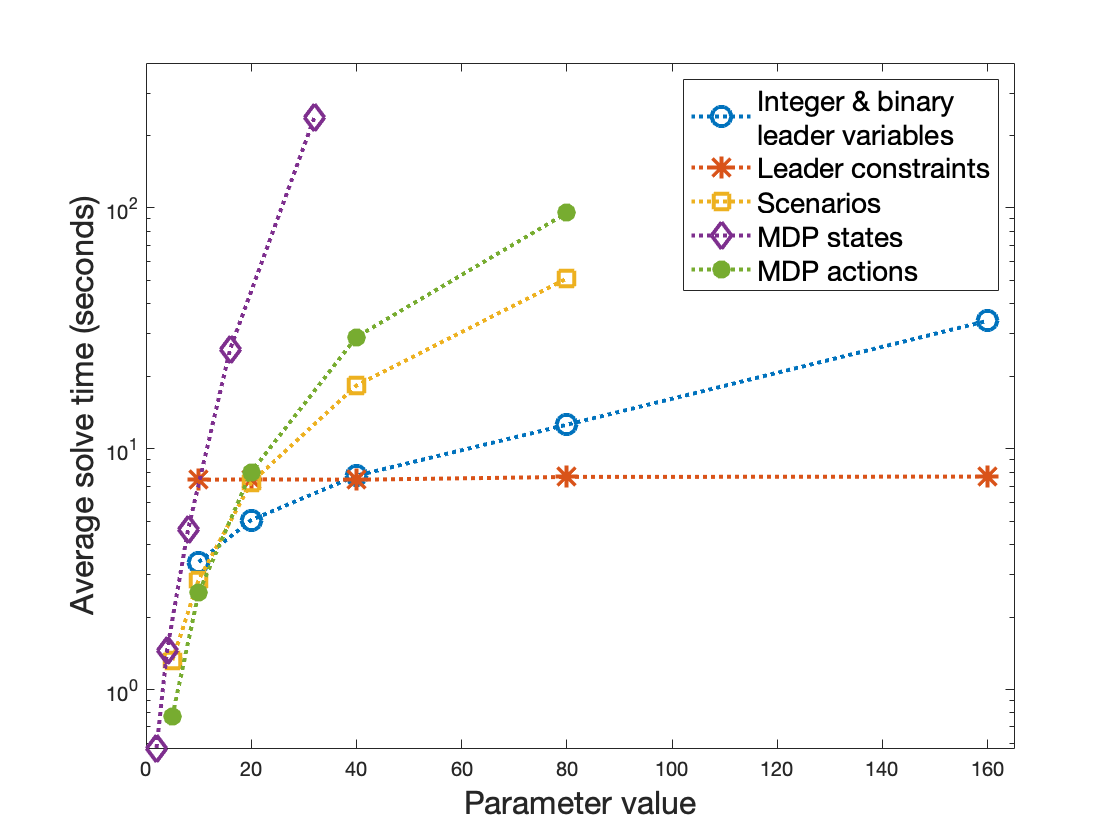} 
    \caption{Trends in average solve time (in seconds) over 100 instances upon varying the number of leader variables, leader constraints, scenarios, MDP states and MDP actions per state.}
    \label{fig:MC-level-C}
\end{figure}

The average solve times for each combination of parameter values are tabulated in Table \ref{table:solvetimes}. Trends in average solve time are illustrated in Figure \ref{fig:MC-level-C}.
The number of leader constraints appears to have little effect on the average solve time. However, average solve time increase with the number of MDP states, number of MDP actions per state, number of scenarios, and number of leader variables. Preliminary tests also showed that varying the mean values used to generate right-hand sides for the leader constraints does not significantly affect average solve times. 

\section{Conclusion}
We presented a modeling framework to integrate strategic and operational decisions made by the same agent. \label{rev:conclusion}\rev{The strength of the proposed framework lies in its generality and applicability to a wide range of application domains.}
Given a static design phase modeled as a MIP and a dynamic operational phase represented by an infinite-horizon MDP, we derived a bilevel optimization formulation that captures the temporal hierarchy of the two decision-phases. The bilevel program, which has a mixed-integer linear leader problem and continuous linear follower problems, can directly be solved using existing computational methods. We provided three examples to illustrate how the modeling approach may be applied in practice, and presented numerical results to illustrate that realistic problem instances can be solved using an existing bilevel solver. 
Future work will explore tailored solution approaches that leverage the MDP structure and can lead to improved computational performance.



\section*{Acknowledgments}
This research was supported by the National Science Foundation grant CMMI 1826323, the National Cancer Institute grant 5R01CA257814-02, and the Office of Naval Research grant N00014-21-1-2262. The authors thank Tyler Perini of US Naval Academy for his helpful comments. 

\bibliographystyle{plainnat} 
\bibliography{SPIMRref}
\end{document}